\documentclass[12pt]{article}

\usepackage[utf8]{inputenc}

\usepackage{amsmath}
\usepackage{amsthm}

\usepackage{amssymb}
\usepackage{amscd}
\usepackage{amsfonts}
\usepackage{amsbsy}

\let\emptyset\varnothing

\usepackage{amssymb,amsmath}

\usepackage{color}
\usepackage{graphicx}

\usepackage{wrapfig}

\usepackage{indentfirst}

\usepackage{makeidx}
\makeindex

\usepackage{ordinalpt}

\usepackage{hyperref}

\newtheorem{definicao}{Definition}[section]

\newtheorem{lema}[definicao]{Lemma}%[section]
\newtheorem{observacao}[definicao]{Remark}%[section]
\newtheorem{teorema}[definicao]{Theorem}%[section]
\newtheorem{corolario}[definicao]{Corollary}%[section]
\newenvironment{prova}{\noindent{\bf Proof:}}{\hfill$\qed$}

\newcommand{\T}{\mathbb{T}}
\newcommand{\R}{\mathbb{R}}
\newcommand{\N}{\mathbb{N}}
\newcommand{\Z}{\mathbb{Z}}
\newcommand{\C}{\mathbb{C}}
\newcommand{\nn}{{\bf n}}
\newcommand{\nm}{{\bf m}}
\newcommand{\nk}{{\bf k}}
\newcommand{\nl}{{\bf l}}
\newcommand{\nj}{{\bf j}}

\newcommand{\np}{{\bf p}}
\newcommand{\nq}{{\bf q}}
\newcommand{\nw}{{\bf w}}
\newcommand{\nx}{{\bf x}}
\newcommand{\nc}{{\bf c}}

\newcommand{\nr}{{\bf r}}
\newcommand{\ny}{{\bf y}}

\newcommand{\sen}{\, {\rm{sin}} \,}

\begin{document}

\title{Optimal Approximation by $sk$-Splines on the Torus }
\author{J. G. Oliveira\footnote{E-mail address: jgaiba@yahoo.com.br} \ and
S. A. Tozoni\footnote{E-mail address: tozoni@ime.unicamp.br}}
\date{}
\maketitle

%%%%%%%%%%%%%%%%%%%%%%%%%%%%%%

\center {\em Instituto de Matemática, Universidade Estadual de Campinas, Rua Sérgio Buarque de Holanda 651,
13083-859, Campinas, SP, Brazil}

%%%%%%%%%%%%%%%%%%%%%%%%%%%%%%
\vspace{.5cm}

{\bf Key words and phrases:} spline, torus, interpolation, approximation, multiplier, Sobolev spaces 

{\bf MSC2010 classification numbers:} 41A15, 41A25, 41A65

\vspace{.5cm}

\begin{abstract}

Fixed a continuous kernel K on the $d$-dimensional torus, we consider a generalization of the univariate $sk$-spline to the torus, associated with the kernel K. It is proved an estimate which provides the rate of convergence of a given function by its interpolating $sk$-splines, in the norm of $L^q$ for functions of the type $f=K*\varphi$ where $\varphi \in L^p$ and $1\leq p \leq 2 \leq q \leq \infty,\ 1/p - 1/q \geq 1/2$. The rate of convergence is obtained for functions f in Sobolev classes and this rate gives optimal error estimate of the same order as best trigonometric approximation, in a special case.
 
\end{abstract}

%||||||||||||||||||||||||||||||||||||||||||||||||||||||||||||||||||||

\section{Introduction}

The $sk$-splines are a natural generalization of the polynomial splines and of the $\cal L$-splines of Micchelli \cite{PM}. The $sk$-splines were introduced and their basic theory developed by A. K. Kushpel. The latest results about convergence of $sk$-splines on the circle in the spaces $L^q$, were obtained by Kushpel in \cite{ku2008-14, ku2008-15}. For an overview of approximation by $sk$-splines see  \cite{PM}.

In Section \ref{sec.resultadosbasicos} we introduce the concept of $sk$-spline on the torus $\T^d$ and we show some basic results. In Section \ref{Fundamental sk-spline} we define the fundamental $sk$-spline, interpolating $sk$-splines and we find conditions for the existence and uniqueness of  interpolating $sk$-splines of a given function. We show that the interpolating $sk$-spline can be obtained from the fundamental $sk$-spline.
In Section \ref{aprox.sksplines} we prove Theorem \ref{teorema002} which provides the rate of convergence of a given function by its interpolating $sk$-splines. The rate of convergence is given in the norm of $L^q(\T^d)$ for functions of the type $f=K\ast \varphi$ where $\varphi \in L^p(\T^d)$ and $K$ is a fixed kernel, for $1\leq p \leq 2 \leq q \leq \infty$ where $p^{-1}-q^{-1}\geq 2^{-1}$. The most important result for our applications is the Corollary \ref{coronucleo}.

Consider fixed a kernel K. Given $\nn=(n_1,\ldots,n_d)\in \N^d$ and $\nk=(k_1,\ldots,k_d)\in \Z^d$ let $\nx_{\nk}=(x_{k_1},\ldots, x_{k_d}),\ x_{k_l}=\pi k_l/n_l$. We denote by $sk_\nn(f,\cdot)$ the unique interpolating $sk$-spline of a function f with set of knots and interpolating points $\Lambda_\nn=\{ \nx_\nk: 0 \leq k_l\leq 2n_l-1, 1\leq l \leq d   \}.  $
In the last section we prove that if $\gamma \in \R,\ \gamma>d$, 
\begin{equation}\label{intnucleo1}
K(\nx)=\sum_{\nl \in \Z^d \setminus \{ {\bf 0} \}}|\nl|^{-\gamma} e^{i\nl\cdot \nx},\ \ \nx \in \T^d,
\end{equation}
where $|\cdot| $ is the norm $|\cdot|_{2}$ or $|\cdot|_{\infty}$ on $\R^d$,  $\nn=(n,\ldots,n)\in\mathbb N^d$, $1\leq p\leq 2\leq q\leq\infty$, with $1/p-1/q\geq 1/2$, then there is a positive constant $C_{p,q}$, independent of $n\in\mathbb N$, such that
\begin{equation}\label{inteq3.1}\sup_{f\in K*U_p}||  f-sk_\nn(f,\cdot) ||_q\leq C_{p,q} n^{-\gamma+d(1/p-1/q)}.\end{equation}
The set $K*U_p=\{K*\varphi:\varphi\in U_p \},$ where $U_p$ is the unit ball of $L^p(\T^d)$, is a Sobolev class on the torus $\T^d$.

It follows from \cite{KST} and \cite{P} that for $1\leq p\leq 2\leq q<\infty$, the $(2n)^d$-width of Kolmogorov of $K*U_p$ verifies
\begin{equation}\label{inteq3.9}
d_{(2n)^d}(K*U_p,L^q) \asymp n^{-\gamma+d(1/p-1/2)}.
\end{equation}

For $\nn=(n,\ldots,n)\in\mathbb N^d$, the dimension of the space $SK(\Lambda_\nn)$ of interpolating $sk$-splines on $\Lambda_\nn$ is $(2n)^d$. 
Comparing \eqref{inteq3.1} and \eqref{inteq3.9} we can see that the rate of convergence by interpolating $sk$-splines is as good as the rate of convergence by subspaces of trigonometric polynomials of the dimension of $SK(\Lambda_\nn)$ on Sobolev classes, when $p=1$ and $q=2$, that is, the rate of convergence is optimal in the sense of $n$-widths.
The construction of the optimal interpolating $sk$-splines is given in Theorem \ref{teo02}.  $SK(\Lambda_\nn)$ is an optimal subspace for the Kolmogorov $(2n)^d$-width of the Sobolev class $K*U_1$ in $L^2$.

In \cite{LK, LK2} the convergence of $sk$-splines for functions in anisotropic Sobolev classes on the torus was studied. These studies were improved in \cite{GKLR2, GKLR}. The best result was obtained in  \cite{GKLR}. It was proved an almost optimal estimate, in the sense of best approximation by trigonometric polynomials, for functions in Sobolev classes, by $sk$-splines, optimal up to a logarithmic factor. In \cite{KGD} it was obtained a similar result  for the case $p=q=1$.

\section{Preliminaries}

If $(a_n)$ and $(b_n)$ are sequences, we write $a_n\gg b_n$ to indicate that there is a constant $C_1>0$ such that $a_n\geq C_1 b_n$ for all $n\in\mathbb N$ and we write $a_n\ll b_n$ to indicate that there is a constant $C_2>0$ such that $a_n\leq C_2 b_n$ for all $n\in\mathbb N$. We write $a_n\asymp b_n$ to indicate that $a_n\ll b_n$ and $a_n\gg b_n$.

The $d-$dimensional torus $\T^d$ is defined as the product of $d$ copies of the quotient group $\R / 2\pi \Z$, or $\T^d = \R / 2\pi \Z \times \R / 2\pi \Z \times \cdots \times \R / 2\pi \Z.$
We can identify $\T^d$ with the $d-$dimensional cube $[-\pi,\pi]^d$ and also with the cartesian product $S^1\times \cdots \times S^1$, of $d$ times the unitary circle $S^1=\{ e^{it}:t \in [-\pi,\pi] \}$. 
We will consider $\T^d$ endowed with the normalized Lebesgue measure $d\nu(\nx)=\left({1}/{(2\pi)^d}\right)dx_1dx_2\cdots dx_d$, where $\left({1}/{2\pi}\right)dt$ is the normalized Lebesgue measure on $S^1$.  

For $\nl=(l_1,\ldots,l_d),\, \nk=(k_1,\ldots,k_d),\, \nj=(j_1,\ldots,j_d)\in \Z^d$ and ${\bf x}=(x_1,\ldots,x_d), \, {\bf y}=(y_1,\ldots,y_d)\in \mathbb R^d$, we denote
 ${\bf x}\cdot {\bf y}=x_1y_1+\cdots+x_dy_d$;
 $\nl\nk=(l_1k_1,\ldots,l_dk_d)$;
 $\nl\equiv \nk\ {\rm mod}(\nj)$ if there is $\np \in \Z^d$ such that $\nl-\nk=\np\nj;$
 ${\bf 0}=(0,0, \ldots, 0)$;
 ${\bf 1}=(1,1, \ldots, 1)$;
 $|{\bf x}|_p=(|x_{1}|^p+|x_2|^p+\cdots+|x_d|^p)^{1/p}\ $ for $ \ 1\leq p < \infty; $
 $|{\bf x}|_{\infty}=\max_{1\leq j \leq d}|x_j|$.

In this paper we consider a arbitrary norm $\nx \rightarrow |\nx|$ on $\R^d$ and we denote by $|\nl|$ the norm of the element $\nl \in \Z^d$.

We denote by $L^p=L^p(\T^d),\ 1\leq p \leq \infty$, the vector space of all mensurable functions $f$ defined on $\T^d$ and with values in $\C$, satisfying 
\begin{eqnarray*}
&& ||f||_p= \left( \int_{\T^d}|f(\nx)|^p d \nu(\nx)\right)^{1/p}< \infty, \ 1 \leq p <\infty,\\
&&||f||_\infty = {\rm ess} \sup_{\nx \in \T^d}|f(\nx)| < \infty. 
\end{eqnarray*}
We write $U_p=\{ f \in L^p(\T^d):||f||_p\leq 1 \}$.

Given $f\in L^1(\T^d)$ we define the Fourier series of the function $f$ by
\begin{equation*}
\sum_{\nm \in \Z^d}\widehat{f}(\nm)e^{i\nm\cdot \nx}.
\end{equation*}
where $$ \widehat{f}(\nm)= \int_{\T^d}f(\nx)e^{-i\nm\cdot \nx}d\nu(\nx).$$

The convolution product of two functions $f$ and $g$ in $L^1(\T^d)$, denoted by $f*g$, is defined by
$$f*g(\nx)=\int_{\T^d}f(\nx-\ny)g(\ny)d\nu(\ny).$$
If $1\leq p,q \leq \infty, \ f \in L^q(\T^d)$ and $g \in L^p(\T^d)$, then the Young Inequality says that $f*g \in L^s(\T^d)$, where $1/s=1/p+1/q-1$, and $$||f*g||_{s}\leq||f||_{q}||g||_{p}.$$

Let $(a_{\bf l})_{{\bf l}\in\mathbb Z^d}$ be a sequence of real numbers such that $a_{\bf l}=a_{{\bf -l}}$ for every ${\bf l}\in\mathbb Z^d$ and \[ \sum_{\bf l\in\mathbb Z^d} |a_{\bf l}|<\infty.\]
Consider the kernel $K({\bf x})$ given by \[K({\bf x})=\sum_{{\bf l}\in\mathbb Z^d}a_{\bf l} e^{i{\bf l}\cdot {\bf x}}.\]
We have that $K$ is a real function, continuous and even. We consider the convolution operator defined for $f\in L^1(\mathbb T^d)$ by \begin{equation*} Tf(x)=K* f(x), \, x\in\mathbb T^d.\end{equation*}
$T$ is a bounded linear operator from $L^p(\mathbb T^d)$ to $L^q(\mathbb T^d)$, for $1 \leq p,q \leq \infty$. For $f\in L^1(\T^d)$ we have
\begin{equation*}\label{1.6}
Tf(\nx)=\sum_{\nl\in \mathbb Z^d} a_{\nl}\hat{f}(\nl)e^{i\nl\cdot\nx},
\end{equation*}
and the norm of T as an operator from $L^p$ to $L^q$ is the norm of T as an operator from  $L^{p'}$ to $L^{q'}$, that is $||T||_{p,q}=||T||_{q',p'}$, 
for every $p,q\in\mathbb R$, $1\leq p,q\leq\infty$, where $p'$ and $q'$ satisfy $1/p+1/p'=1/q+1/q'=1$. We denote $$K*U_p=\{K*f: f \in U_p \}.$$

For $l, N \in \N$ we define
$$A_l= \{ \nk \in \Z^d:|\nk|_2\leq l \},\ \ A_l^*= \{ \nk \in \Z^d:|\nk|_{\infty}\leq l \},$$
$${\cal H}_l=[e^{i\nk\cdot \nx}: \nk \in A_l\setminus A_{l-1}],\ {\cal H}_l^*=[e^{i\nk\cdot \nx}: \nk \in A_l^*\setminus A^*_{l-1}],$$
where $A_{-1}=A_{-1}^*=\emptyset$ and $[f_j: j\in \varGamma ]$ denotes the vector space generated by the functions $f_j:\T^d \longrightarrow \C$, with $j$ in the set of indexes $\varGamma$. We denote by ${\cal H}$ the vector space generated by the family $\{ e^{i\nk\cdot \nx}: \nk \in \Z^d \}$ which is dense in $L^p(\T^d)$ for $1 \leq p < \infty$.
Then by \cite{FR} and \cite{MI}, there are positive constants $C_1$, $C_2$ and $C_3$ satisfying 
$$\frac{2 \pi^{d/2}}{\varGamma(d/2)}l^{d-1}-C_2l^{d-2}\leq \dim{{\cal H}_l} \leq \frac{2 \pi^{d/2}}{\varGamma(d/2)}l^{d-1}+C_1l^{d-2}. $$
It is easily to verify that there is a positive constant $C$ such that for every $l,N\in\mathbb N$,
$$\dim{{\cal H}_l}^*=(2l+1)^d-(2(l-1)+1)^d \asymp l^{d-1}.$$
In particular $\dim{{\cal H}_l}\asymp \dim{{\cal H}_l^*} \asymp l^{d-1}$.

Consider a function $\lambda: [0,\infty) \rightarrow \R$ and let $\nx \rightarrow |\nx|$ be a norm on $\R^d$. For each $\nk \in \Z^d$ we define $\lambda_\nk=\lambda(|\nk|)$.
We denote by $\Lambda$ the linear operator defined for $\varphi \in {\cal H}$ by 
$$\Lambda \varphi = \sum_{\nk \in \Z^d} \lambda_\nk \widehat {\varphi}(\nk)e^{i\nk\cdot \nx}.$$

Let $\Lambda=\{ \lambda_\nk \}_{\nk \in \Z^d}$, $\lambda_\nk \in \C$, and $1\leq p,q \leq \infty$. If for any $\varphi \in L^p(\T^d)$ there is a function $f=\Lambda\varphi \in L^q(\T^d)$ with formal Fourier  expansion given by
$$f \thicksim \sum_{\nk \in \Z^d}\lambda_\nk \widehat \varphi(\nk)e^{i\nk\cdot \nx}$$
such that $||\Lambda||_{p,q}= \sup \{ ||\Lambda\varphi||_{q}: \varphi \in U_p \}< \infty$, we say that $\Lambda$ is a bounded multiplier operator from $L^p(\T^d)$ into $L^q(\T^d)$, with norm $||\Lambda||_{p,q}$.

%&&&&&&&&&&&&&&&&&&&&&&&&&&&&&&&&&&&&&&&&&&&&&&&&&&&&&&&&&&&&&&&&&&&&&&&&&&&&&&&&&&&&&&&&&&&&&&&&
%&&&&&&&&&&&&&&&&&&&&&&&&&&&&&&&&&&&&&&&&&&&&&&&&&&&&&&&&&&&&&&&&&&&&&&&&&&&&&&&&&&&&&&&&&&&&&&&&
%&&&&&&&&&&&&&&&&&&&&&&&&&&&&&&&&&&&&&&&&&&&&&&&&&&&&&&&&&&&&&&&&&&&&&&&&&&&&&&&&&&&&&&&&&&&&&&&&&&&&&&&&&&&&&&&&&&&&&&&&&&&&&&&
\section{Basic results}
\label{sec.resultadosbasicos}

Let $\nn=(n_1, \ldots, n_d)\in \mathbb{N}^d$. For $\nk =(k_1, \ldots, k_d) \in \Z^d$ we denote $x_{k_{l}}={\pi k_l}/{n_l}$,\ $ 1 \leq l \leq d$ and
$\mathbf{x_{k}}=(x_{k_1}, \ldots, x_{k_d})$. We also denote 
$$\Omega_{\nn}=\{ \mathbf{j}=(j_1, \ldots, j_d) \in \mathbb{Z}^d: \ 0\leq j_l \leq 2n_l-1,\  1\leq l \leq d \},$$
$$\Lambda_{\nn}=\{\mathbf{x_k} : \mathbf{k} \in \Omega_{\mathbf{n}} \},\ \ \ \ \   N= \# \Omega_{\nn} = \# \Lambda_{\nn}= 2^dn_1n_2\cdots n_d.$$

The real vector space of all continuous functions $f:\T^d\rightarrow \R$ endowed with the norm of the uniform convergence will be denoted by $C(\T^d)$.

For a fixed kernel $K \in C(\mathbb{T}^d)$, a $sk$-spline on $\Lambda_{\nn}$ is a function represented in the form $$sk_\nn(\mathbf{x})=c+\sum_{\mathbf{k} \in \Omega_{\mathbf{n}}}c_\mathbf{k} K(\mathbf{x}-\mathbf{x_k}),$$
where the coefficients $c$, $c_\mathbf{k} \in \mathbb{R},\ \mathbf{k} \in \Omega_\mathbf{n}$, satisfy the condition $$\sum_{\nk \in \Omega_\nn} c_\nk= 0.$$ The points $\mathbf{x_k}$ are the knots of the $sk$-spline $sk_\nn(\mathbf{x})$.

The real vector space of all $sk$-splines on $\Lambda_{\nn}$, associated with the kernel $K$ will be denoted by $SK(\Lambda_{\nn})$. As the vector space $V$ generated by the set of functions $\{1,K(\nx-\nx_{\nk}),\ \nk \in \Omega_\nn \}$ has dimension at most $N+1$ and $SK(\Lambda_\nn)$ is the subspace of $V$ formed by the functions whose coefficients satisfy the condition $\sum_{\nk \in \Omega_\nn}c_\nk=0$, then ${\rm{dim}}$ $SK(\Lambda_\nn) \leq N$.

The next four lemmas will not be proved, because they are of simple verification.

\begin{lema}\label{lema01}
Let $\mathbf{l} \in \mathbb{Z}^d$. Then
\[ \sum_{{\mathbf{k}}\in \Omega_{\nn}} e^{i \mathbf{l \cdot x_k}} =\left\{
\begin{array}{cc}
N,& \mathbf{l} \equiv {\bf 0}\ {\rm mod}(2\mathbf{n}),\\
0,& {\textrm{{\rm otherwise},}}
\end{array}
\right.\ \ \ and\ \ \ \int_{\T^d} e^{i \mathbf{l \cdot x}}d \nu (\nx) =\left\{
\begin{array}{cc}
0,& \mathbf{l} \neq {\bf 0},\\
1,& \nl={\bf 0}.
\end{array}
\right.\]
\end{lema}

\begin{lema}\label{coro01}
For every $\mathbf{l} \in \mathbb{Z}^d$,
\[ \sum_{{\mathbf{\nk}}\in \Omega_{\nn}} {\cos} (\mathbf{l\cdot x_{k}}) =\left\{
\begin{array}{cc}
N,& \mathbf{l} \equiv \mathbf{0}\ {\rm mod}(2\mathbf{n}),\\
0,& \textrm{{\rm otherwise,}}
\end{array}
\right.\ \ \ {{and}} \ \ \ \sum_{{\mathbf{k}}\in \Omega_{\nn}} {\sen}(\mathbf{l\cdot x_{k}}) =0.\]
\end{lema}

\begin{lema}\label{coro02}
For every $\mathbf{l,j} \in \mathbb{Z}^d$ we have that
\[ \sum_{{\mathbf{k}}\in \Omega_{\nn}} ({\cos}(\mathbf{j\cdot x_{k}}))({\cos}(\mathbf{l\cdot x_{k}})) =\left\{
\begin{array}{cc}
N,& \mathbf{l+j} \equiv {\bf 0}\ {\rm mod}(2\mathbf{n})\ {\rm and}\\
 & \mathbf{l-j} \equiv {\bf 0}\ {\rm mod}(2\mathbf{n}),\ \\
N/2,& \mathbf{l+j} \equiv {\bf 0}\ {\rm mod}(2\mathbf{n})\ {\rm or}\\
 & \mathbf{l-j} \equiv {\bf 0}\ {\rm mod}(2\mathbf{n}),\ \ \ \\
0,& \textrm{{\rm otherwise,}}
\end{array}
\right.\]

\begin{equation*}
\sum_{{\mathbf{k}}\in \Omega_{\nn}} ({\sen}(\mathbf{j\cdot x_{k}}))({\sen}(\mathbf{l\cdot x_{k}})) =\left\{
\begin{array}{cc}
N/2,& \mathbf{l+j} \not \equiv {\bf 0}\ {\rm mod}(2\mathbf{n})\ {\rm and}\\
 & \mathbf{l-j} \equiv {\bf 0}\ {\rm mod}(2\mathbf{n}),\ \\
-N/2,& \mathbf{l+j} \equiv {\bf 0}\ {\rm mod}(2\mathbf{n})\ {\rm and}\\
 & \mathbf{l-j} \not\equiv {\bf 0}\ {\rm mod}(2\mathbf{n}),\ \ \\
0,& \textrm{{\rm otherwise,}}
\end{array}
\right.
\end{equation*}
and
\[\sum_{{\mathbf{k}}\in \Omega_{\nn}} ({\cos}(\mathbf{j\cdot x_{k}}))({\sen}(\mathbf{l\cdot x_{k}})) =0.\]
\end{lema}

\begin{definicao}
For $K\in C(\T^d)$, $\mathbf{j}\in \Z^d$ and $\mathbf{x} \in \mathbb{T}^d$, we define
$$\lambda_{\mathbf{j}}(\mathbf{x})=\sum_{\mathbf{k}\in \Omega_\nn}e^{i\mathbf{j\cdot  x_k}}K(\mathbf{x-x_{k}}),$$
$$\rho_\nj(\nx)=\frac{2}{N}{\rm Re}(\lambda_\nj(\nx))=\frac{2}{N}\sum_{\mathbf{k}\in \Omega_\nn}(\cos{(\mathbf{j\cdot  x_k}}))K(\mathbf{x-x_{k}})$$
and
$$\sigma_\nj(\nx)=\frac{2}{N}{\rm Im}(\lambda_\nj(\nx))=\frac{2}{N}\sum_{\mathbf{k}\in \Omega_\nn}(\sen{(\mathbf{j\cdot  x_k}}))K(\mathbf{x-x_{k}}).$$
\end{definicao}
\begin{lema}\label{lema3}
Let $\np,\nj \in \Z^d$. Then for every $\nx \in \mathbb{T}^d$,
$$\rho_{2\nn \np+\nj}(\nx)=\rho_\nj(\nx),\ \ \ 
\rho_{2\nn \np -\nj}(\nx)=\rho_\nj(\nx),\ \ \ 
\rho_{-\nj}(\nx)=\rho_\nj(\nx),$$
$$\sigma_{2\nn \np+\nj}(\nx)=\sigma_\nj(\nx),\ \ \ 
\sigma_{2\nn \np -\nj}(\nx)=-\sigma_\nj(\nx),\ \ \ 
\sigma_{-\nj}(\nx)=-\sigma_\nj(\nx).$$
\end{lema}

\begin{teorema}\label{lema05}
Consider a kernel $K$ given by
$K(\mathbf{x})=\sum_{\mathbf{l}\in \Z^d}a_{\nl}e^{i\mathbf{l\cdot  x}},$
where $\left( a_\nl \right)_{\nl \in \Z^d}$ is a sequence of real numbers such that
$\sum_{\mathbf{l}\in \Z^d}|a_{\nl}|<\infty$
and $a_\nl=a_{-\nl}$ for every ${\mathbf{l}\in \Z^d}$. Then $\ K$ is a real function, continuous, even and for every  ${\mathbf{j}\in \Z^d}$,
$$\rho_{\mathbf{j}}(\mathbf{x})=\sum_{{\mathbf{p}}\in \Z^d} (a_{2\mathbf{ np+j}}{\cos}\,((2\nn\np+\nj)\cdot \nx)+a_{{2\mathbf{ np-j}}}{\cos}((2 \nn\np-\nj)\cdot \nx)),$$
$$\sigma_{\mathbf{j}}(\mathbf{x})=\sum_{{\mathbf{p}}\in \Z^d} (a_{2\mathbf{ np+j}}{\sen}\,((2\nn\np+\nj)\cdot \nx)-a_{{2\mathbf{ np-j}}}{\sen}((2 \nn\np-\nj)\cdot \nx)).$$
\end{teorema}
\begin{prova}
For each $m\in \N $, let
$K_m(\mathbf{x})=\sum_{|\mathbf{l}|_2\leq m}a_{\nl}e^{i\mathbf{l\cdot  x}}.$
We have that $(K_m)_{m\in\N}$ is a Cauchy sequence in $C(\T^d)$ and as $C(\T^d)$ is complete, there is a function $K\in C(\T^d)$ such that $K_m \rightarrow K$ uniformly. 
We have that
\begin{eqnarray}\label{defk}
K(\mathbf{x})
&=&\sum_{\mathbf{l}\in \Z^d}a_{\nl}{\cos}({\mathbf{l\cdot  x}})
\end{eqnarray}
and $K$ is a real and even function. Fix $\nj \in \Z ^d$ and let  
$$A_{\nj}=\{ \nl \in \Z^d : \nl+\nj\equiv {\bf 0}\ {\rm mod}(2\nn)\}=\{ 2\nn \mathbf{p}-\nj:\mathbf{p} \in \Z^d \},$$
$$B_{\nj}=\{ \nl \in \Z^d : \nl-\nj\equiv {\bf 0}\ {\rm mod}(2\nn)\}=\{ 2\nn \mathbf{p}+\nj:\mathbf{p} \in \Z^d \}.$$
From Lemma \ref{coro02} we have that
\begin{equation}\label{0.14}
\sum_{{\mathbf{k}}\in \Omega_{\nn}} ({\cos}(\mathbf{j\cdot x_{k}}))({\cos}(\mathbf{l\cdot x_{k}})) =\left\{
\begin{array}{cc}
N,& \nl \in A_\nj \cap B_\nj, \\
N/2,& \nl \in A_\nj \triangle B_\nj,\\
0,& \nl \in (A_\nj \cup B_\nj)^{c}
\end{array}
\right.
\end{equation}
\begin{equation}\label{0.014}
\sum_{{\mathbf{k}}\in \Omega_{\nn}} ({\sen}(\mathbf{j\cdot x_{k}}))({\sen}(\mathbf{l\cdot x_{k}})) =\left\{
\begin{array}{cc}
N/2,& \nl \in B_\nj\backslash A_\nj, \\
-N/2,& \nl \in A_\nj\backslash B_\nj,\\
0,& \nl \in (A_\nj \triangle B_{\nj})^{c}.
\end{array}
\right.
\end{equation}
Now using (\ref{defk}) we obtain
\begin{eqnarray}\label{0.15}
\lambda_{\mathbf{j}}(\mathbf{x})&=&\sum_{\mathbf{k}\in \Omega_\nn}e^{i\mathbf{j\cdot  x_k}}K(\mathbf{x-x_{k}})\nonumber\\
&=&\sum_{\nl \in \Z ^d}a_\nl \sum_{\mathbf{k}\in \Omega_\nn}({{\cos}(\mathbf{j\cdot  x_k}}))({\cos}(\mathbf{l}\cdot(\mathbf{x-x_{k}})))\nonumber\\
&+&i\sum_{\nl \in \Z ^d}a_\nl \sum_{\mathbf{k}\in \Omega_\nn}({{\sen}(\mathbf{j\cdot  x_k}}))({\cos}(\mathbf{l}\cdot(\mathbf{x-x_{k}}))).
\end{eqnarray}
Thus by \eqref{0.15}, from Lemma \ref{coro02} and from \eqref{0.14}
\begin{eqnarray*}
\rho_\nj(\mathbf{x})&=&\frac{2}{N} \sum_{\nl \in \Z ^d}a_\nl \sum_{\mathbf{k}\in \Omega_\nn}({{\cos}(\mathbf{j\cdot  x_k}}))({\cos}(\mathbf{l}\cdot\mathbf{x-l\cdot x_{k}}))\\
&=&\frac{2}{N} \sum_{\nl \in \Z ^d}a_\nl({\cos}(\mathbf{l}\cdot\mathbf{x})) \sum_{\mathbf{k}\in \Omega_\nn}({{\cos}(\mathbf{j\cdot  x_k}}))({\cos}(\mathbf{l}\cdot\mathbf{x_k}))\\
&=&\frac{2}{N}\sum_{\nl \in A_\nj\cap B_\nj}Na_\nl {\cos}(\nl \cdot \nx)+\frac{2}{N}\sum_{\nl \in A_\nj\triangle B_\nj}\frac{N}{2}a_\nl {\cos}(\nl \cdot \nx)\\
&=&\sum_{2\nn\np+\nj \in A_\nj\cap B_\nj, \np \in \Z^d}a_{2\nn\np+\nj} {\cos}((2\nn\np+\nj) \cdot \nx)\\
&+&\sum_{2\nn\np-\nj \in A_\nj\cap B_\nj, \np \in \Z^d}a_{2\nn\np-\nj} {\cos}((2\nn\np-\nj) \cdot \nx)\\
&+&\sum_{2\nn\np+\nj \in A_\nj\triangle B_\nj, \np \in \Z^d}a_{2\nn\np+\nj} {\cos}((2\nn\np+\nj) \cdot \nx)\\
&+&\sum_{2\nn\np-\nj \in A_\nj\triangle B_\nj, \np \in \Z^d}a_{2\nn\np-\nj} {\cos}((2\nn\np-\nj) \cdot \nx)\\
&=&\sum_{2\nn\np+\nj \in A_\nj\cup B_\nj, \np \in \Z^d}a_{2\nn\np+\nj} {\cos}((2\nn\np+\nj) \cdot \nx)\\
&+&\sum_{2\nn\np-\nj \in A_\nj\cup B_\nj, \np \in \Z^d}a_{2\nn\np-\nj} {\cos}((2\nn\np-\nj) \cdot \nx)\\
&=&\sum_{\np \in \Z^d}a_{2\nn\np+\nj} {\cos}((2\nn\np+\nj) \cdot \nx)+\sum_{ \np \in \Z^d}a_{2\nn\np-\nj} {\cos}((2\nn\np-\nj) \cdot \nx).
\end{eqnarray*}
In an analogous way, using \eqref{0.15}, Lemma \ref{coro02} and \eqref{0.014} we obtain
\begin{eqnarray*}
\sigma_\nj(\nx)
&=&\sum_{\np \in \Z^d}a_{2\nn\np+\nj} {\sen}((2\nn\np+\nj) \cdot \nx)-\sum_{ \np \in \Z^d}a_{2\nn\np-\nj} {\sen}((2\nn\np-\nj) \cdot \nx),
\end{eqnarray*}
and this concludes the proof.
\end{prova}

\section{Fundamental sk-spline}
\label{Fundamental sk-spline}
\begin{definicao}\label{def01}
Suppose that $\rho_\nj({\bf 0})\neq 0$ for all $\nj \in \Omega_\nn,$  $\nj \neq {\bf 0}$. We define $\widetilde{sk}_\nn$ by
$$\widetilde{sk}_\nn(\mathbf{x})=\frac{1}{N}+ \frac{1}{N}\sum_{\mathbf{j} \in \Omega_{\nn}^{*}}\frac{\rho_\mathbf{j}(\mathbf{x})}{\rho_\mathbf{j}(\mathbf{0})}$$
where $\Omega_\nn^*=\Omega_\nn\setminus \{(0,\ldots,0) \}. $
\end{definicao}
\begin{lema}\label{lema4}
The function $\widetilde{sk}_\nn$ is a $sk$-spline.
\end{lema}
\begin{prova}
We have, by the definition of $\rho_{\nj}(\nx)$ that
\begin{eqnarray*}
\widetilde{sk}_\nn(\mathbf{x})&=&\frac{1}{N}+ \frac{1}{N}\sum_{\mathbf{j} \in \Omega_{\nn}^{*}}\frac{\rho_\mathbf{j}(\mathbf{x})}{\rho_\mathbf{j}(\mathbf{0})}\\
&=&\frac{1}{N}+\sum_{\mathbf{k} \in \Omega_{\nn}}\left( \frac{2}{N^2} \sum_{\nj \in \Omega_\nn^{*}}\frac{1}{\rho_\mathbf{j}(\mathbf{0})} ({\cos}(\nj \cdot \nx_\nk))\right)  K(\nx-\nx_\nk)\\
&=&\frac{1}{N}+\sum_{\mathbf{k} \in \Omega_{\nn}}c_\nk  K(\nx-\nx_\nk)
\end{eqnarray*}
and by Lemma \ref{coro01}
$$\displaystyle \sum_{\nk \in \Omega_\nn} c_\nk = \frac{2}{N^2}\sum_{\nj \in \Omega_\nn^{*}}\frac{1}{\rho_\nj({\bf 0})}\sum_{\nk \in \Omega_\nn} \cos(\nj \cdot \nx_\nk)=0.$$
Thus $\widetilde{sk}_\nn$ is a $sk-$spline by definition.
\end{prova}

The sk-spline $\widetilde{sk}_\nn$ will be called fundamental sk-spline.
\begin{lema}\label{lema6}
If $\rho_\nj({\bf{0}})\neq 0$ for all $\nj \in \Omega_\nn^{*}$, then the $sk$-spline $\widetilde{sk}_\nn$ satisfy
$$
\widetilde{sk}_\nn(\nx_\nk)=\left\{
\begin{array}{cc}
1, & \nk={\bf{0}}, \\ 
0,& \nk \in \Omega_\nn ^{*}.\\
\end{array}
\right.
$$
\end{lema}
\begin{prova}
From Theorem \ref{lema05} we have that
\begin{eqnarray*}
\rho_\nj(\nx_\nl)&=&\sum_{\np \in \Z^d}\left( a_{2\nn\np+\nj}{\cos}((2\nn\np+\nj) \cdot \nx_\nl)+a_{2\nn\np-\nj} {\cos}((2\nn\np-\nj) \cdot \nx_\nl)\right) \\
&=&\sum_{\np \in \Z^d}\left( a_{2\nn\np+\nj}{\cos}(\nj \cdot \nx_\nl)+a_{2\nn\np-\nj} {\cos}(\nj \cdot \nx_\nl)\right) \\
&=&({\cos}(\nj \cdot \nx_\nl))\sum_{\np \in \Z^d}\left( a_{2\nn\np+\nj}+a_{2\nn\np-\nj} \right)\\
&=&({\cos}(\nj \cdot \nx_\nl))\rho_\nj({\bf{0}}),\\
\end{eqnarray*}
then by Lemma \ref{coro01}
$$\begin{array}{lcl}
\widetilde{sk}_\nn(\mathbf{x_k})&=&\displaystyle{\frac{1}{N}+\frac{1}{N}\sum_{\mathbf{j} \in \Omega_\nn^{*}}\frac{({\cos}(\nj \cdot \nx_\nk))\rho_\mathbf{j}(\mathbf{0})}{\rho_\mathbf{j}(\mathbf{0})}}
=\displaystyle{\frac{1}{N}\sum_{\mathbf{j} \in \Omega_\nn}{\cos}(\nk \cdot \nx_\nj)}
= \displaystyle{\left\{\begin{array}{cc}
1,  & \nk ={\bf{0}},\\
0,& \nk\in \Omega_\nn^{*},
\end{array}
\right.}\\
\end{array}$$
and then we proved the lemma.
\end{prova}
\begin{definicao}
Let $f$ be a function defined on $\T^d$ and let $\{\ny_\nj: \nj\in \Omega_\nn \} \subset \T^d$. If there are constants $c^*, c_\nk^{*}\in \mathbb{R},$ such that
$${sk}_\nn(f,\mathbf{y}_\nj)=c^*+\sum_{\mathbf{k}\in \Omega_\nn}c_{\nk}^{*}K(\mathbf{\mathbf{y}_j-x_k})=f(\mathbf{\mathbf{y}_\nj}),\ \  \nj \in \Omega_\nn,$$
we say that the $sk$-spline
$${sk}_\nn(f,\nx)=c^*+\sum_{\mathbf{k}\in \Omega_\nn}c_{\nk}^{*}K(\mathbf{\mathbf{x}-x_k})$$
is an interpolating $sk$-spline of $f$ with knots $\nx_\nk$ and interpolation points $\ny_\nk$.
\end{definicao}

\begin{teorema}\label{teo02}
Suppose $\rho_\mathbf{j}(\mathbf{0})\neq 0$ for any $\mathbf{j}\in \Omega_\nn^*$. Then for any function f defined on $\T^d$, there is an unique interpolating $sk$-spline of $f$ with knots and interpolation points $\nx_\nk, \nk\in \Omega_\nn$, that can be written in the form
$$sk_\nn(f,\mathbf{x})=\sum_{\mathbf{k}\in \Omega_\nn}f(\mathbf{x_k})\widetilde{sk}_\nn(\mathbf{x}-\mathbf{x_k}).$$
\end{teorema}
\begin{prova}
Let $c_\nk,\  \nk \in \Omega_\nn$ be the coefficients of the $sk$-spline $\widetilde {sk}_\nn$ that were obtained in the proof of Lemma \ref{lema4} and let
\begin{eqnarray*}
\displaystyle sk_\nn(f,\nx)&=& \sum_{\nk \in \Omega_\nn} \frac{f(\nx_\nk)}{N} + \sum_{\nl \in \Omega_\nn} \left( \sum_{\nk \in \Omega_\nn} c_\nl f(\nx_\nk)\right)K(\nx-\nx_\nl)
= d+\sum_{\nl \in \Omega_\nn} d_\nl K(\nx-\nx_\nl). 
\end{eqnarray*}
Since $\sum_{\nl \in \Omega_\nn} c_\nl=0$, it follows that $\sum_{\nl \in \Omega_\nn} d_\nl=0$, and thus $sk_\nn(f,\cdot)$ is a $sk$-spline.
Applying Lemma \ref{lema6} we obtain that
$$sk_\nn(f,\mathbf{x_l})=\sum_{\mathbf{k}\in \Omega_\nn}f(\mathbf{x_k})\widetilde{sk}_\nn(\mathbf{x_l-x_k})=f(\mathbf{x_l})$$
for any $\mathbf{l}\in \Omega_\nn$. Then we can conclude that $sk_\nn(f,\cdot)$ is an interpolating $sk$-spline of $f$ with knots and interpolation points $\nx_\nk.$

Let $\{\mathbf{w}_j :1\leq j \leq N \}$ be an enumeration of $\Lambda_\nn$. Then for every function $f $ on $\T^d$, there are constants $c_1, c_2, \ldots, c_{N+1} \in \R$ satisfying $\sum_{l=1}^{N}c_l=0$, such that
$$sk_\nn(f,\nx)=c_{N+1}+\sum_{l=1}^{N}{c}_l K(\nx-{\nw}_l)$$
is an interpolating $sk$-spline of $f$. Given $y_1, y_2, \ldots, y_N \in \R,$ let
$$\alpha_j=\left( \prod_{1\leq l \leq N, \ \l\neq j} |\nw_j-\nw_l|_2 \right)^{-1},\ 1\leq j \leq N,$$
and $g:\T^d \rightarrow \R$ defined by 
$$g(\nx)= \sum_{j=1}^{N}y_j\alpha_j \prod_{1\leq l \leq N \ \l\neq j} |\nx-\nw_l|_2,\ \ \nx \in \T^d .$$
Thus $g(\nw_k)=y_k$, for $1\leq k \leq N$. Then there are $c_1, \ldots, c_N,c_{N+1}\in \R$ such that the $sk-$spline
$$sk_{\nn}(\nx)=c_{N+1}+ \sum_{l=1}^{N}c_lK(\nx-\nw_l)$$
is an interpolating $sk-$spline of $g$, that is,
$$c_{N+1}+ \sum_{l=1}^{N}c_lK(\nw_k-\nw_l)=g(\nw_k)=y_k,\ \ 1\leq k\leq N.$$

Let
$$\mathbf{\widetilde{K}}=
 \left( \begin{array}{cccc}
 K(\nw_1-\nw_1)& \cdots &K(\nw_1-\nw_N)& 1\\
 \vdots & \vdots & \vdots & \vdots \\
 K(\nw_N-\nw_1)& \cdots &K(\nw_N-\nw_N) & 1\\
 1& \cdots& 1 & 0  \\
 \end{array} \right),$$
${\bf{C}}=(c_1, c_2, \ldots, c_{N+1}),\ {\bf {Y}}=(y_1, y_2, \ldots, y_{N+1}) \in \R^{N+1}$ and

$$\mathbf{\widetilde{K}}=
 \left( 
 \begin{array}{cc}
 \mathbf{K}&\mathbf{u} \\
 \mathbf{u}^t& 0
 \end{array} \right),\ \  \mathbf{u}=\left(\begin{array}{c}
 1\\
 \vdots \\
 1
 \end{array}\right),\ \ \mathbf{C}^t=\left(\begin{array}{c}
 c_1\\
 \vdots \\
 c_{N+1}
 \end{array}\right),\ \ \mathbf{Y}^t=\left(\begin{array}{c}
 y_1\\
 \vdots \\
 y_{N+1}
 \end{array}\right), $$
where $\mathbf{u}$ is a $N\times 1$ matrix.
Let $W=\left\{ {\bf{C}}=(c_1, c_2, \ldots, c_{N+1}) \in \R^{N+1}: \sum_{l=1}^{N}{c}_l=0\right\}.$
Then for every ${\bf {Y}}\in \R^{N+1}$ with $y_{N+1}=0$, there is ${\bf {C}}\in W$ such that ${\bf {\widetilde K}}{\bf {C}}^{t}={\bf {Y}}^{t}$.

Now we consider the linear map $T:W \rightarrow \R^{N+1}$ defined by
$T({\bf {C}})=({\bf {\widetilde K}}{\bf {C}}^{t})^t.$
We can conclude that $T$ is injective, since $T$ is linear and $\dim W=N=\dim $ Im(T).

Let $f$ be a function on $\T^d$ and suppose that there are $C=(c_1, c_2, \ldots, c_{N+1})$, $\overline C=(\bar c_1, \bar c_2, \ldots, \bar c_{N+1})\in W$ such that
$${sk}_\nn(f,\nx)= c_{N+1}+\sum_{l=1}^{N}{c}_l K(\nx-\mathbf{w}_l)\ \ \ and\ \ \ \overline{sk}_\nn(f,\nx)=\bar{c}_{N+1}+\sum_{l=1}^{N}\bar{\nc_l} K(\nx-\mathbf{w}_l)$$
are two interpolating $sk$-splines of $f$. If $F=(f(\nw_1), \ldots, f(\nw_N),0)$, then we have
$T(C)=({\bf {\widetilde K}}{\bf {C}}^{t})^t=F$ and
$T(\overline{C})=({\bf {\widetilde K}}{\bf {\overline C}}^{t})^t=F.$
Since $T$ is injective, it follows that $C=\overline C$, that is, ${sk}_\nn(f,\nx)=\overline{sk}_\nn(f,\nx)$, for all $ \nx \in \T^d$.
\end{prova}

\begin{observacao}\label{obs227}
Let $K$ be a kernel satisfying the conditions of the Theorem \ref{lema05} and $\nn=(n_1,n_2,\ldots,n_d)\in \N^d$. Suppose $\rho_\nj(0)\neq 0$ for all $\nj\in \Omega_\nn$. Then the vector space  $SK(\Lambda_\nn)$ of all $sk$-splines on $\Lambda_\nn$ and associated with the kernel $K$ has dimension $N=2^dn_1n_2\ldots n_d$. In particular if $n_1=n_2=\cdots=n_d=n$ we have $\dim(SK(\Lambda_\nn))=(2n)^d$.
\end{observacao}

\section{Approximation by sk-splines}
\label{aprox.sksplines}

In this section we will prove the main result of this paper, the Theorem \ref{teorema002}. This theorem says how a function of the type $f=K*\phi$, for $\phi \in L^p(\T^d)$, can be approximated by the $sk$-splines $sk_\nn(f,\cdot)$ in the space $L^q(\T^d)$, where $1\leq p\leq 2 \leq q \leq \infty$ with $1/p -1/q \geq 1/2$. But for our applications, the most interesting result is the Corollary \ref{coronucleo}, since its hypothesis can be easily verified.

In all results of this section, we consider a kernel $K$ as in the Theorem \ref{lema05} and such that $\rho_\nj(0)\neq 0$ for all $\nn \in \N^d$ and $\nj \in \Omega_\nn$.

The following result can be easily verified.

\begin{lema}\label{lemaaux}
For $\nj \in \Omega_\nn$ and $\nl \in \Z^d$ we have that
$\rho_\nl(\nx-\nx_\nj)=\rho_\nl(\nx)\cos (\nl\cdot\nx_\nj)+\sigma_\nl(\nx)\sen (\nl\cdot\nx_\nj).$
\end{lema}

\begin{lema}\label{lemaaux2}
For every $\nl \in \mathbb{Z}^d$, $\nl\not \equiv {\bf 0}\ {\textrm {\rm mod}}(2\nn)$ and $\nx \in \T^d$,
$$\sum_{\nj \in \Omega_\nn}e^{i\nl\cdot \nx_\nj}\widetilde{sk}_\nn(\nx-\nx_\nj)=\frac{\lambda_\nl(\nx)}{\rho_{\nl}({\bf0})}.$$
\end{lema}

\begin{prova}
Firstly we will prove the result for the real part. 
% Colocar I aqui
Consider the sets $A_\nl$ and $B_\nl$ introduced in the proof of Theorem \ref{lema05}.

Using Lemma \ref{lemaaux} together with Lemmas \ref{coro01} and \ref{coro02}, the equation (\ref{0.14}) and the fact that $\nl\not \equiv {\bf 0}\ {\textrm {\rm mod}}(2\nn)$, we have that
\begin{eqnarray}\label{eq120}
\sum_{\nj\in \Omega_\nn}(\cos(\nl\cdot\nx_\nj))\widetilde{sk}_\nn(\nx-\nx_\nj)
&=&\sum_{\nj\in \Omega_\nn}(\cos(\nl\cdot\nx_\nj))\left\{\frac{1}{N}+ \frac{1}{N}\sum_{\mathbf{k} \in \Omega_{\nn}^{*}}\frac{\rho_\mathbf{k}(\mathbf{x-x_j})}{\rho_\mathbf{k}(\mathbf{0})} \right\}\nonumber\\
&=&\frac{1}{N}\sum_{\mathbf{k} \in \Omega_{\nn}^{*}}\frac{\rho_\mathbf{k}(\mathbf{x})} {\rho_\mathbf{k}(\mathbf{0})}\sum_{\nj\in \Omega_\nn}(\cos(\nl\cdot\nx_\nj))(\cos(\nk\cdot\nx_\nj))  \nonumber \\
&=&\frac{1}{N}\sum_{\mathbf{k} \in \Omega_{\nn}^{*} \cap(A_\nl \cap B_\nl)}\frac{\rho_\mathbf{k}(\mathbf{x})} {\rho_\mathbf{k}(\mathbf{0})}\sum_{\nj\in \Omega_\nn}(\cos(\nl\cdot\nx_\nj))(\cos(\nk\cdot\nx_\nj))\nonumber\\
&+& \frac{1}{N}\sum_{\mathbf{k} \in \Omega_{\nn}^{*} \cap(A_\nl \Delta B_\nl)}\frac{\rho_\mathbf{k}(\mathbf{x})} {\rho_\mathbf{k}(\mathbf{0})}\sum_{\nj\in \Omega_\nn}(\cos(\nl\cdot\nx_\nj))(\cos(\nk\cdot\nx_\nj))\nonumber\\
&=& \sum_{\mathbf{k} \in \Omega_{\nn}^{*} \cap(A_\nl \cap B_\nl)}\frac{\rho_\mathbf{k}(\mathbf{x})} {\rho_\mathbf{k}(\mathbf{0})}
+ \frac{1}{2}\sum_{\mathbf{k} \in \Omega_{\nn}^{*} \cap(A_\nl \Delta B_\nl)}\frac{\rho_\mathbf{k}(\mathbf{x})} {\rho_\mathbf{k}(\mathbf{0})}.
\end{eqnarray}

If $\nk \in B_\nl$ then there is $\np \in \Z^d$ such that $\nk=2\nn\np+\nl$ and thus $\rho_\nk(\nx)=\rho_{2\nn\np+\nl}(\nx)=\rho_\nl(\nx)$ by Lemma \ref{lema3}. In an analogous way, if $\nk \in A_\nl$ we can conclude that $\rho_\nk(\nx)=\rho_\nl(\nx)$. Then, by (\ref{eq120}) we have that
\begin{equation}\label{eq1.6t}
\sum_{\nj\in \Omega_\nn}(\cos(\nl\cdot\nx_\nj))\widetilde{sk}_\nn(\nx-\nx_\nj)=\frac{\rho_\nl(\nx)}{\rho_{\nl}({\bf0})}\left(\#(\Omega_{\nn}^{*} \cap(A_\nl \cap B_\nl))+\frac12 \#(\Omega_{\nn}^{*} \cap(A_\nl \Delta B_\nl)) \right). 
\end{equation}

Let $\nl=(l_1,\ldots,l_d) \in \Z^d$. For each $1\leq j \leq d$, there is an unique $q_j$ and an unique $r_j$ satisfying $q_j, r_j \in \Z, 0\leq r_j \leq 2n_j-1$ and $l_j=2 n_j q_j + r_j$. Then $ \nl=2\nn \nq + \nr $ 
where $\nq =(q_1, \ldots ,q_d) \in \Z^d$ and $\nr=(r_1, \ldots ,r_d) \in \Omega_{\nn}$, so
$$B_\nl=\{ 2\nn \np + \nl: \np \in \Z^d  \}=\{ 2\nn (\np+\nq) + \nr: \np \in \Z^d  \}=\{ 2\nn \np + \nr: \np \in \Z^d  \}=B_\nr.$$

In an analogous way we have $A_\nl=\{ 2\nn \np - \nr: \np \in \Z^d  \}=A_\nr.$
As $\rho_\nl(\nx)=\rho_\nr(\nx)$ by Lemma \ref{lema3}, we obtain by (\ref{eq1.6t}) that
$$\sum_{\nj\in \Omega_\nn}(\cos(\nl\cdot\nx_\nj))\widetilde{sk}_\nn(\nx-\nx_\nj)=\frac{\rho_\nr(\nx)}{\rho_{\nr}({\bf0})}\left( \#\left(\Omega_\nn^{*} \cap ( A_{\nr} \cap B_{\nr}) \right) + \frac12 \#(\Omega_\nn^{*} \cap ( A_{\nr} \Delta B_{\nr})  ) \right).$$
Then it is enough to prove the result for $\nl \in \Omega_\nn^*$.

Let $\nl=(l_1, \ldots, l_d), \nk=(k_1, \ldots, k_d) \in \Omega_\nn^*$. Then $\nl-\nk \equiv{\bf 0}\mod(2\nn)$ if and only if $\nk=\nl$, and $\nl+\nk \equiv {\bf 0} \mod(2\nn)$ if and only if $k_j=l_j=0$ and $k_j=2n_j-l_j$ if $l_j\neq 0$. Then
$\nk \in \Omega_\nn^{*} \cap ( A_{\nl} \cap B_{\nl}) $ if and only if $\nk=\nl$ and $l_j\in \{ 0,n_j\}$ for all $1\leq j \leq d$; $\nk \in \Omega_\nn^{*} \cap ( B_{\nl} \setminus A_{\nl})$ if and only if $\nk=\nl$ and $l_j\not\in \{ 0,n_j\}$ for some $1\leq j \leq d$; $\nk \in \Omega_\nn^{*} \cap ( A_{\nl} \setminus B_{\nl}) $ if and only if $k_j=l_j=0,\ k_j=2n_j-l_j$ if $l_j\neq 0$ and $l_j\not \in \{ 0,n_j\}$ for some $1\leq j \leq d$. Let
$${\cal{A}}=\{ \nl \in \Omega_\nn^{*}: l_j \in \{ 0,n_j\}\ {\rm{for \ all}}\ 1\leq j \leq d \},$$
$${\cal B}=\{ \nl \in \Omega_\nn^{*}: l_j \not\in \{ 0,n_j\}\ {\rm{for \ all}}\ 1\leq j \leq d \}.$$
Then $\Omega_{\nn}^{*}={\cal{A}} \cup {\cal{B}},\  {\cal{A}}\cap {\cal{B}}=\emptyset \ $  and  
$$\ \# (\Omega_{\nn}^{*} \cap ( A_{\nl} \cap B_{\nl}) )=\left\{
\begin{array}{cc}
1,& \nl \in {\cal{A}},\\
0,& \nl \in {\cal{B}},
\end{array}
\right. ,\ \ \ \ \ \  \# (\Omega_{\nn}^{*} \cap ( A_{\nl} \Delta B_{\nl}) )=\left\{
\begin{array}{cc}
0,& \nl \in {\cal{A}},\\
2,& \nl \in {\cal{B}}.
\end{array}
\right.$$
Then it follow from (\ref{eq1.6t}) that for all $\nl \in \Omega_\nn^*$ we have
\begin{equation}\label{eq1.7t}
\sum_{\nj\in \Omega_\nn}(\cos(\nl\cdot\nx_\nj))\widetilde{sk}_\nn(\nx-\nx_\nj)=\frac{\rho_\nl(\nx)}{\rho_{\nl}({\bf0})}.
\end{equation}

In an analogous way, for the imaginary part, we obtain that for all $\nl \in \Z^d$, $\nl\neq {\bf 0}$,
$$\sum_{\nj\in \Omega_\nn}(\sen(\nl\cdot\nx_\nj))\widetilde{sk}_\nn(\nx-\nx_\nj)=\frac{\sigma_\mathbf{\nl}(\mathbf{x})} {\rho_\mathbf{\nl}(\mathbf{0})},$$
and this concludes the proof.
\end{prova}

\begin{observacao}\label{obs2.2.3}
Let $|\cdot|$ be a norm on $\mathbb R^d$ and let $K\in C(\T^d)$ be a kernel as in Theorem \ref{lema05}, such that  $a_{\nl}=a_{\nk}$ if $\nl,\nk\in\mathbb Z^d$ and $|\nl|=|\nk|$. Given $\nk=(k_1,\ldots,k_d)$, $\np=(p_1,\ldots,p_d)$, ${\bf i}=(i_1,\ldots,i_d)\in\mathbb Z^d$, let
$\overline\nk=((-1)^{i_1}k_1,\ldots,(-1)^{i_d}k_d)$ and $\overline\np=((-1)^{i_1}p_1,\ldots,(-1)^{i_d}p_d)$.
Then $|2\nn\np+\overline\nk|=|2\nn\overline\np+\nk|$
and so $a_{2\nn\np+\overline\nk}=a_{2\nn\overline\np+\nk}$.

Given $\nj=(j_1,\ldots,j_d)\in \left(\N \cup \{ 0 \} \right)^d =\{0,1,2,\ldots\}^d$, let
\[D_{\nj}= \left\{\np=(p_1,\ldots,p_d)\in\mathbb Z^d:\, |p_i|=j_i, \ i=1,2,\ldots,d \right\}.\]
Then
\[
\begin{array}{lcl}
\displaystyle \sum_{\np\in\mathbb Z^d} a_{2\nn\np+\overline\nk}&=&\displaystyle \sum_{\nj\in\left(\N \cup \{ 0 \} \right)^d}\sum_{\np\in D_{\nj}} a_{2\nn\np+\overline\nk}
=\displaystyle \sum_{\nj\in\left(\N \cup \{ 0 \} \right)^d}\sum_{\np\in D_{\nj}}a_{2\nn\overline\np+\nk}\\
&=&\displaystyle \sum_{\nj\in\left(\N \cup \{ 0 \} \right)^d}\sum_{\np\in D_{\nj}}a_{2\nn\np+\nk}
=\displaystyle \sum_{\np\in\mathbb Z^d} a_{2\nn\np+\nk}.
\end{array}
\]
\end{observacao}

\begin{observacao}\label{obs2.3.4}
Consider a kernel $K$ given by $K(\nx)=\sum_{\nl \in \Z^d}a_{\nl}e^{i\nl \cdot \nx },$
such that $a_{\nl}\geq 0$, for all $\nl \in \Z^d$ and $\sum_{\np \in \Z^d}a_{2\nn\np-\nk}\leq Ca_{2\nn-\nk},$ for all $\nn=(n_1,n_2,\ldots,n_d)\in \N^d$ and every $\nk=(k_1,\ldots,k_d)\in\mathbb Z^d$, with $0\leq k_j \leq n_j$, for $j=1,2,\ldots,d$, where $C$ is a positive constant independent of $\nn$ and $\nk$. Then
$\sum_{\nl \in \Z^d}a_{\nl}<\infty.$
If $a_\nl=a_{-\nl}$ for all $\nl\in\mathbb Z^d$, then by Theorem \ref{lema05}, the kernel $K$ is a real, continuous and even function.
\end{observacao}

\begin{lema}\label{lemaaux19}
Let $|\cdot |$ be a norm on $\mathbb R^d$ and let $K$ be the kernel given by 
$K(\nx)=\sum_{\nl\in\mathbb Z^d}a_{\nl}e^{i\nl\cdot\nx}$,
where $(a_{\nl})_{\nl\in\mathbb Z^d}$ is a sequence with  $a_{\nl}=a_{\nk}$ if $|\nl|=|\nk|$ and $a_{\nl}\geq a_{\nk}>0$ if $|\nl|\geq|\nk|$, for $\nl,\nk\in\mathbb Z^d$. 
Suppose that there is a positive constant $C$ such that for every $\nn\in\mathbb N^d$ and all $\nk=(k_1,\ldots,k_d) \in \mathbb{Z}^d$, with $0\leq k_j \leq n_j$ for $j=1,2,\ldots, d$, we have
$\sum_{\np\in \Z^d}a_{2\nn\np-\nk}\leq Ca_{2\nn-\nk}.$
Let $$\theta_{\nn,\nl}(\nx)=\displaystyle e^{i\nl\cdot \nx}-\sum_{\nj \in \Omega_\nn} e^{i \nl\cdot \nx_\nj} \widetilde{sk}_{\nn}\left(\nx-\nx_\nj \right).\vspace{0.3cm}
$$
Then for $ \nl \in \mathbb{Z}^d$, $\tilde{\nl}=(|l_1|, \ldots, |l_d|),$
\[\begin{array}{lcl}
|\theta_{\nn,\nl}(\nx)| &\leq&\left\{
\begin{array}{lcl}
\displaystyle 4C\frac{a_{2\nn-\tilde{\nl}}}{a_\nl},&&{ 0}< | \nl |\leq |\nn|,\\
|e^{i\nl\cdot \nx}- 1|,&& {\rm{for\ }}\ \nl \equiv {\bf 0}\ {\rm mod}\ (2\mathbf{n}),\\
4,&& {\rm{for\ all}}\ \nl.
\end{array}
\right.
\end{array}
\]
\end{lema}
\begin{prova}
Let $\mu_{\nn,\nl}(\nx)$ be the real part of $\theta_{\nn,\nl}(\nx)$. For $\nl \equiv {\bf 0}\ {\rm mod}\ (2\mathbf{n})$, by the Definition \ref{def01}, Lemmas \ref{coro01} and \ref{lemaaux} we have
\[ 
\mu_{\nn,\nl}(\nx)=\displaystyle\cos(\nl\cdot \nx)-\sum_{\nj \in \Omega_\nn} (\cos\left(\nl\cdot \nx_\nj \right))\widetilde{sk}_\nn\left(\nx-\nx_\nj\right) = \displaystyle \cos(\nl\cdot \nx)- 1.
\]

For $\nl \not\equiv {\bf 0}\ {\rm mod}\ (2\mathbf{n})$, using Lemma \ref{lemaaux2} we have
\[
\begin{array}{lcl}
\vspace{0.3cm}\mu_{\nn,\nl}(\nx)&=&\displaystyle \cos(\nl\cdot \nx)-\frac{\rho_\nl(\nx)}{\rho_\nl({\bf{0}})}=\displaystyle \frac{\rho_\nl({\bf 0})\cos(\nl\cdot \nx)-\rho_{\nl}(\nx)}{\rho_\nl({\bf{0}})}.
\end{array}
\]
% onde \[s(\nx)=\sum_{\nm \in \Z^d} \{ a_{2\nm\nn-\nl}\cos((2\nm\nn-\nl)\cdot\nx)+a_{2\nm\nn+\nl}\cos((2\nm\nn+\nl)\cdot\nx) \}.\]
Thus, as in the Theorem \ref{lema05} we have $\rho_\nl(\nx)\leq \rho_\nl({\bf{0}})$ for all $\nx$ and for all $\nl\in \Z^d$, then 
\[
\begin{array}{lcl}
|\mu_{\nn,\nl}(\nx)|&=&\displaystyle\left| \displaystyle \frac{\rho_\nl({\bf 0})\cos(\nl\cdot \nx)-\rho_{\nl}(\nx)}{\rho_\nl({\bf{0}})}\right|
\leq\displaystyle \frac{\rho_\nl({\bf 0})+\rho_{\nl}(\nx)}{\rho_\nl({\bf{0}})} \leq \frac{2\rho_\nl({\bf 0})}{\rho_\nl({\bf{0}})} =2.
\end{array}
\]
Suppose now that ${ 0}< | \nl |\leq |\nn|$ and let $\tilde{\nl}=(|l_1|, \ldots, |l_d|)$. Thus using the hypothesis and the Remark $\ref {obs2.2.3}$, we obtain
\[
\begin{array}{lcl}
|\mu_{\nn,\nl}(\nx)|&=& \displaystyle\left| \displaystyle \frac{\rho_\nl({\bf 0})\cos(\nl\cdot \nx)-\rho_{\nl}(\nx)}{\rho_\nl({\bf{0}})}\right|
\leq \displaystyle \frac{2\rho_\nl({\bf{0}})}{2a_\nl}
\leq 2C \displaystyle \frac{a_{2\nn-\tilde{\nl}}}{a_\nl}\\
\end{array}
\]

The imaginary part of $\theta_{\nn,\nl}(\nx)$ is given by  
\[ \phi_{\nn,\nl}(\nx)=\sen(\nl\cdot \nx)-\sum_{\nj \in \Omega_\nn} (\sen\left(\nl\cdot \nx_\nj \right))\widetilde{sk}_\nn\left(\nx-\nx_\nj\right). \]
The estimate for the imaginary part is analogous to the real part. Considering the estimates obtained for $\mu_{\nn,\nl}(\nx)$ and $\phi_{\nn,\nl}(\nx)$, we obtain the desired estimate for $\theta_{\nn,\nl}(\nx)$.
\end{prova}

\begin{lema}\label{lemaaux20}
Let $K$ be a kernel as in Lemma \ref{lemaaux19}. Then for each $1\leq p<\infty$, there is a positive constant $C$, depending only on $p$, such that
\[\sum_{\nl \in\mathbb Z^d} a_\nl^p|\theta_{\nn,\nl}(\nx)|^p\leq C\sum_{|\nl|\geq |\nn|} a_\nl^p.\]
\end{lema}
\begin{prova} Since $a_\nl\geq a_\nk>0$ if $|\nk|\geq |\nl|,\ \nk,\nl\in \Z^d$, using Lemma \ref{lemaaux19} and taking $\tilde{\nl}=(|l_1|,\ldots, |l_d|)$ we have
\[
\begin{array}{lcl}
\displaystyle\sum_{\nl \in\mathbb Z^d} a_\nl^p|\theta_{\nn,\nl}(\nx)|^p
&\leq&\displaystyle \sum_{{0}< |\nl| \leq |\nn|} a_{\nl}^p |\theta_{\nn,\nl}(\nx)|^p+\sum_{|\nl|\geq |\nn|}a_{\nl}^p|\theta_{\nn,\nl}(\nx)|^p\vspace{0.3cm}\\
&\leq&\displaystyle \sum_{{ 0}< |\nl| \leq |\nn|}a_{\nl}^p 4^pC^p\left( \dfrac{a_{2\nn-\tilde{\nl}}}{a_{\nl}}\right)^p+\sum_{|\nl | \geq |\nn|}a_{\nl}^p4^p \vspace{0.3cm}\\
&=&\displaystyle 4^pC^p  \sum_{{0}< |\nl| \leq |\nn|} a_{2\nn-\tilde{\nl}}^p+ 4^p\sum_{|\nl|\geq |\nn|}a_{\nl}^p. \vspace{0.3cm} \\
\end{array}
\]
For each $\nl \in \Z^d, \nl \neq {\bf 0}$, let
$D_\nl=\{\nk=(k_1,\ldots,k_d)\in \Z^d: |k_j|=|l_j|, \ 1\leq j \leq d \}.$
Then
\[
\begin{array}{lcl}
\displaystyle\sum_{\nl \in\mathbb Z^d} a_\nl^p|\theta_{\nn,\nl}(\nx)|^p
&\leq&\displaystyle 4^pC^p \sum_{0<|\nl| \leq |\nn|} \sum_{\nk\in D_\nl} a_{2\nn-{\nk}}^p+ 4^p\sum_{|\nl|\geq |\nn|}a_{\nl}^p \vspace{0.3cm} \\
&\leq&\displaystyle 4^pC^p2^d \sum_{{0}< |\nl| \leq |\nn|} a_{2\nn-{\nl}}^p+ 4^p\sum_{|\nl|\geq |\nn|}a_{\nl}^p \vspace{0.3cm} \\
&\leq&\displaystyle C_1\left( \sum_{|\nn| \leq |\nj| \leq 3|\nn|}a_{\nj}^p+\sum_{|\nl|\geq |\nn|} a_{\nl}^p \right)\leq 2C_1\sum_{|\nl|\geq |\nn|}a_{\nl}^p,
\end{array}
\]
completing the proof of the lemma.
\end{prova}

\begin{teorema} \label{teorema002}
Let $|\cdot|$ be a norm on $\mathbb R^d$ and let $K$ be a kernel given by
\[K(\nx)=\sum_{\nl \in \Z^d}a_\nl e^{i\nl\cdot \nx},\]
where $(a_{\nl})_{\nl\in\mathbb Z^d}$ is a sequence that satisfies $a_\nl=a_{\nk}$ if $|\nl|=|\nk|$ and
$a_\nl\geq a_\nk>0$ if $|\nk|\geq |\nl|$, for $\nk,\nl\in \Z^d$. Suppose that
\[\sum_{\np\in \Z^d}a_{2\nn\np-\nk}\leq Ca_{2\nn-\nk},\]
for all $\nn\in\mathbb N^d$ and all $\nk=(k_1,\ldots,k_d)$ with $0 \leq k_j \leq n_j$ for $j=1, 2,\ldots, d,$ where $C$ is a positive constant that is independent of $\nn$ and $\nk$. Then for $1\leq p\leq 2\leq q\leq \infty$, with $p^{-1}-q^{-1}\geq 2^{-1}$, we have
\[ \displaystyle  \sup_{f\in K*U_p}|| f-{ sk_\nn}(f,\cdot)||_q\leq C\left(\sum_{|\nl|\geq |\nn|}a_\nl^{qp(q-p)^{-1}}\right)^{p^{-1}-q^{-1}}. \]
\end{teorema}
\begin{prova}
Let $p\in \R,\ 1\leq p \leq 2$ and let $p'$  such that $1/p+{1}/{p'}=1$. Given $f\in K * U_p$, $\phi\in U_p$ such that $f=K*\phi$, by Theorem \ref{teo02},
$$\begin{array}{lcl}
\displaystyle \sigma_{\nn}(f,\nx)&=&f(\mathbf{x})-sk_\nn(f,\mathbf{x})\\
&=&\displaystyle \int_{\mathbb{T}^d}\left(K(\mathbf{x-y})-\sum_{\mathbf{k} \in \Omega_\nn}K(\mathbf{x_k-y})\widetilde{sk}_\nn(\mathbf{x-x_k})\right)\phi(\mathbf{y})d\nu(\mathbf{y})\\
&=&\displaystyle \int_{\mathbb T^d}  \Phi_{\nn}(\nx,\ny)\phi(\ny)d\nu(\ny),
\end{array}$$
where $\Phi_{\nn}(\nx,\ny)=K(\nx,\ny)-\sum_{\nk\in\Omega_{\nn}} K({\nx}_{\nk}-\ny)\widetilde{sk}_{\nn}(\nx-\nx_{\nk})$. Thus by Hölder inequality we have that
\begin{equation}\label{eq101}
\displaystyle \left| f(\nx)-{ sk}_\nn(f,\nx) \right| 
\leq \displaystyle ||\phi||_p  || \Phi_{\nn}(\nx,\cdot)||_{p'}.
\end{equation}
Since $1\leq p\leq 2$, it follows from Hausdorff-Young inequality that
$$\left| \left| \Phi_{\nn}(\nx,\cdot) \right|\right|_{p'}\leq \left( \sum_{\nl \in \Z^d}|b_\nl|^p \right)^{1/p}, $$
where for $\nl \in \Z^d$,
$b_\nl =  \int_{\T^d} \Phi_{\nn}(\nx,\ny)e^{-i\nl \cdot \ny}  d\nu(\ny) .$
By Lemma \ref{lema01} we have that
$$ \begin{array}{lcl}
\displaystyle \int_{\mathbb{T}^d}K(\mathbf{x-y})e^{-i \nl \cdot \ny}d\ny
&=& \displaystyle  \sum_{\nj \in \Z^d }a_\nj e^{i\nj \cdot \nx} \int_{\mathbb{T}^d}e^{-i (\nl+\nj) \cdot \ny}d\nu(\ny)
= a_\nl e^{-i\nl \cdot \nx}
\end{array}$$
and in an analogous way
$$ \begin{array}{lcl}
\displaystyle \int_{\mathbb{T}^d}\left( \sum_{\nk \in \Omega_\nn }K(\mathbf{x_\nk-y})\widetilde{sk}_\nn(\nx-\nx_\nk)e^{-i \nl \cdot \ny}\right)d\nu(\ny)
&=&\displaystyle \sum_{\nk \in \Omega_\nn } a_\nl e^{-i\nl \cdot \nx_\nk}\widetilde{sk}_\nn(\nx-\nx_\nk).
\end{array}$$
Thus
\begin{eqnarray*}
 b_\nl &=&a_\nl \left(e^{-i\nl\cdot \nx}- \sum_{\nk \in \Omega_\nn} e^{-i\nl\cdot \nx_\nk} \widetilde{sk}_\nn(\nx-\nx_\nk) \right)
 =a_\nl  \theta_{\nn,-\nl}(\nx)=a_{-\nl}  \theta_{\nn,-\nl}(\nx).
\end{eqnarray*}
Using Lemma \ref{lemaaux20} we obtain
\begin{eqnarray}\label{eq102}
 \left| \left| \Phi_{\nn}(\nx,\cdot) \right| \right|_{p'}
&\leq & \left( \sum_{\nl \in \Z^d} a_{-\nl}^p \left| \theta_{\nn,-\nl}(\nx)\right|  ^p \right)^{1/p}
\leq  C\left( \sum_{|\nl|\geq |\nn|}a_\nl^p \right) ^{1/p}.
\end{eqnarray}
For $\phi \in L^p(\T^d)$ we define
$$ T\phi (\nx)=\int_{\T^d}\Phi_{\nn}(\nx,\ny) \phi(\ny)d\nu(\ny). $$
By inequalities (\ref{eq101}) and (\ref{eq102}) we conclude that $T$ is a bounded operator from $L^p(\T^d)$ to $L^{\infty}(\T^d)$ and that
\begin{equation}\label{eqI}
||T||_{p,\infty}\leq C\left( \sum_{|\nl|\geq |\nn|}a_\nl^p \right)^{1/p}.
\end{equation}
By duality, $T$ is bounded from $L^1(\T^d)$ to $L^{p'}(\T^d)$ and
\begin{equation}\label{eqII}
||T||_{1,p'}\leq C\left( \sum_{|\nl|\geq |\nn|}a_\nl^p \right)^{1/p}.
\end{equation}
Applying the Riesz-Thorin Interpolation Theorem we have $1\leq (p_t^{-1}-q_t^{-1})^{-1}\leq p$ and 
%tomando $r=p_t$ e $s=q_t$, temos $1\leq r \leq 2 \leq s \leq \infty$ e $r^{-1}-s^{-1}\geq 2^{-1}$. Portanto $T$ é limitado de $L^r(\T^d)$ em $ L^s(T^d) $ e
$$||T||_{p_t,q_t}\leq C\left( \sum_{|\nl|\geq |\nn|}a_\nl^{q_t p_t(q_t-p_t)^{-1}} \right)^{p^{-1}_{t}-q^{-1}_{t}}. $$

If $1 \leq r \leq 2$, $2 \leq s \leq \infty$ and ${1}/{r}-{1}/{s}\geq {1}/{2}$, then there are $0\leq t \leq 1$ and $1 \leq p \leq 2$ such that ${1}/{r}=1-t+{t}/{p}$ and $ 
{1}/{s}=({1-t})/{p'},$  
that is, $r=p_t$ and $s=q_t$. 
\end{prova}

\begin{lema}\label{lema2.3.7}
Let $a:[0,+\infty)\rightarrow\mathbb R$ be a decreasing and positive function and $|\cdot| = |\cdot|_p$ for some $\ 1 \leq  p \leq \infty$. For each $\np \in \Z^d$, let $a_{\np}=a(|\np|)$.
Suppose that there is a constant $c_1>0$ such that for each $\nn\in\mathbb N^d$, 
\begin{equation}\label{eq3.6}
\sum_{\np\in\Z^d}a_{2\nn\np}\leq c_1a_{2\nn}. \ \ \ 
\end{equation}
Then there is a constant $c_2>0$ such that for each $\nn\in\mathbb N^d$ and $\nk \in \Z^d$ with $|\nk|\leq |\nn|$, we have
\begin{equation}\label{eq3.7}
\sum_{\np\in\Z^d}a_{2\nn\np-\nk}\leq c_2a_{2\nn-\nk}. \ \ \ 
\end{equation}

\end{lema}

\begin{prova}
Fix $\nn=(n_1,\ldots,n_d)\in\mathbb N^d$. By Remark \ref{obs2.2.3} it is enough to consider $\nk=(k_1,\ldots,k_d)\in \Z^d$ with $|\nk|\leq|\nn|$ and $0\leq k_j \leq n_j$, for each $j=1,2,\ldots,d.$
Let $\np=(p_1,\ldots,p_d)\in\Z^d$. 
For each $1\leq j\leq d$, if
\[\psi_j(p_j)=\tilde p_j=\left\{
\begin{array}{lcl}
p_{j}-1&,&p_j>0,\\
p_j&,&p_j\leq 0,
\end{array}
\right.\]
we define 
$\psi(\np)=\tilde \np=(\psi_1(p_1),\ldots,\psi_d(p_d)),$
and $\psi$ is well defined as a function from $\mathbb Z^d$ to $\mathbb Z^d$.
As
$\begin{array}{lcl}
|2{n_j p_j-k_j}| \geq |2n_j\tilde p_j|
\end{array}$
we have
$$
|2\nn\np-\nk|_p
\geq  \left( |2n_1 \tilde p_1|^p+\cdots+|2n_d\tilde p_d|^p \right)^{1/p}  
=|2\nn\tilde \np|_p,\ 1\leq p <\infty,$$
$$
|2\nn\np-\nk|_\infty
\geq \max\{   |2 n_1 \tilde p_1|,\ldots,|2n_d\tilde p_d| \}
=|2\nn\tilde \np|_\infty.$$
Thus
$|2\nn\np-\nk|\geq|2\nn\tilde \np|$
and consequently
$
a_{2\nn\np-\nk}\leq a_{2\nn\tilde \np}.
$
Since the cardinality of $\psi^{-1}(\{\nk\})$ is at most $2^d$ for all $\nk\in\mathbb Z^d$, by \eqref{eq3.6} 
$$
\sum_{\np\in\mathbb Z^d}a_{2\nn\np-\nk}
\leq\sum_{\np\in\mathbb Z^d}a_{2\nn\tilde{\np}}
\leq2^d\sum_{\tilde{\np}\in\mathbb Z^d}a_{2\nn\tilde{\np}}
\leq2^d c_1 a_{2\nn}
\leq2^d c_1 a_{2\nn-\nk},
$$
and this concludes the proof.
\end{prova}

The next result is consequence of Theorem \ref{teorema002} and Lemma \ref{lema2.3.7}.

\begin{corolario}\label{coronucleo}
Let $a:[0,+\infty)\rightarrow\mathbb R$ be a decreasing and positive function and $|\cdot| = |\cdot|_p$ for some $\ 1 \leq  p \leq \infty$. For each $\np \in \Z^d$ let $a_{\np}=a(|\np|)$.
Consider the kernel $K$ given by
\[K(\nx)=\sum_{\nl \in \Z^d}a_\nl e^{i\nl\cdot \nx},\]
such that
$$\sum_{\np\in \Z^d}a_{2\nn\np}\leq Ca_{2\nn},$$
where $C$ is a positive constant independent of $\nn\in\mathbb N^d$. Then there is a positive constant $\overline{C}$, such that for each $1\leq p\leq 2\leq q\leq \infty$, with $p^{-1}-q^{-1}\geq 2^{-1}$ and all $\nn\in\mathbb N^d$, we have
\[ \displaystyle  \sup_{f\in K*U_p}|| f-{ sk_\nn}(f,\cdot)||_q\leq \overline{C}\left(\sum_{|\nl|\geq |\nn|}a_\nl^{qp(q-p)^{-1}}\right)^{p^{-1}-q^{-1}}. \]
\end{corolario}

%######################################################################################################################################################################
%######################################################################################################################################################################
%###################################################################################################

%%%%%%%%%%%%%%%%%%%%%%%%%%%%%%%%%%%%%%%%%%%%%%%%%%%%%%%%%%%%%%%%%%%%%
\section{Approximation of finitely differentiable functions}

\begin{teorema}\label{teo311}
For $\gamma \in \R,\ \gamma>d$, let
\[K(\nx)=\sum_{\nl \in \Z^d \setminus \{ {\bf 0} \}}|\nl|^{-\gamma} e^{i\nl\cdot \nx},\ \ \nx \in \T^d,\]
where $|\cdot|=|\cdot|_{2}$ or $|\cdot|=|\cdot|_{\infty}$. 
For $n\in\mathbb N$, let $\nn=(n,\ldots,n)\in\mathbb N^d$. Then, for $1\leq p\leq 2\leq q\leq\infty$, with $1/p-1/q\geq 1/2$, there is a positive constant $C_{p,q}$, independent of $n\in\mathbb N$, such that

\begin{equation}\label{eq3.1}\sup_{f\in K*U_p}||  f-sk_\nn(f,\cdot) ||_q\leq C_{p,q} n^{-\gamma+d(1/p-1/q)}.\end{equation} 

\end{teorema}
\begin{prova}
Let $\alpha \in \R$, $\alpha>0$. Using the function $f(x)=(x-1)^\alpha/x^\alpha$, $x\geq 2$ we obtain 
\begin{eqnarray}\label{auxda35}
(j-1)^{-\alpha}\leq 2^{\alpha}j^{-\alpha}, \ j\geq 2.
\end{eqnarray}
Fix $\nn=(n, n, \ldots, n)$ and for each $j\in \N$ let $B_j=\{ \nl \in \mathbb{Z}^d : j-1\leq |\nl |< j \}$. Then $\Z ^d= \bigcup_{j=1}^{\infty}B_j$. If $\np \in B_j$, then $j-1\leq|\np|< j$ and thus $j^{-\gamma}< |\np|^{-\gamma}\leq (j-1)^{-\gamma}$. Let $a_{\nl}=|\nl|^{-\gamma}$ for $\nl \in \Z^d\setminus \{ {\bf 0} \}$ and $a_{{\bf 0}}=0$.
We have $\dim{{\cal H}_l}\asymp \dim{{\cal H}_l^*} \asymp l^{d-1}$ and then the cardinality of $B_j$ satisfies
\begin{equation}\label{eqcontagem}
\# B_j \leq C j^{d-1},\ j \in \N,
\end{equation}
where $C$ is a positive constant independent of $j$. Since
$2\nn\np=2n\np,$
by (\ref{auxda35}) and (\ref{eqcontagem})
\begin{eqnarray*}
\sum_{\np \in \Z^d}a_{2\nn\np}&=&\sum_{j=2}^{\infty}\sum_{\np \in B_j}|2n\np|^{-\gamma}
\leq\sum_{j=2}^{\infty}\sum_{\np \in B_j}(2n)^{-\gamma}(j-1)^{-\gamma}\\
&\leq&C\sum_{j=2}^{\infty}(2n)^{-\gamma}j^{d-1}(j-1)^{-\gamma}
\leq2^{\gamma}C(2n)^{-\gamma}\sum_{j=2}^{\infty}j^{d-1-\gamma}.
\end{eqnarray*}
Since $\gamma>d$, then $d-1-\gamma<0$ and
$\sum_{j=2}^{\infty}j^{d-1-\gamma}\leq \int_{1}^{\infty}t^{d-1-\gamma}dt$. Thus, since $d{-\gamma}<0$,
\begin{eqnarray}\label{umaestrela}
\sum_{\np \in \Z^d}a_{2\nn\np}
&\leq& 2^{\gamma}C(2n)^{-\gamma}\lim_{m \rightarrow \infty }\int_{1}^{m}t^{d-1-\gamma}dt
= \frac{2^{\gamma}C |{\bf 1}|^\gamma}{\gamma-d}a_{2\nn}=C_1 a_{2\nn}.
\end{eqnarray}
\\
Therefore the hypothesis of Corollary \ref{coronucleo} is satisfied.

Let $r=p^{-1}-q^{-1}$ and $s=r^{-1}$. Then using (\ref{auxda35}) and ({\ref{eqcontagem}})
\begin{eqnarray*}
 \sum_{|\nl|\geq|\nn|} (a_{\nl})^{s}
&\leq& \sum_{j=|\nn|+1}^\infty\sum_{\nl \in B_j} (j-1)^{-s\gamma} 
\leq C\sum_{j=|\nn|}^\infty (j+1)^{d-1}j^{-s\gamma} 
\leq 2^{d-1}C\sum_{j=|\nn|}^\infty j^{d-1-s\gamma}. \\
\end{eqnarray*}
We have $1\leq p \leq 2$ and thus $r=1/p-1/q \leq 1$ and $s\geq 1$. Then $d-1-s\gamma<0$ and therefore $j^{d-1-s\gamma}\leq\int_{j-1}^{j}t^{d-1-s\gamma}dt$. We obtain that
\begin{eqnarray*}
\sum_{|\nl|\geq|\nn|} (a_{\nl})^{s}&\leq&  2^{d-1}C\sum_{j=|\nn|}^\infty \int_{j-1}^{j}t^{d-1-s\gamma}dt
=-2^{d-1}C\frac{(|\nn|-1)^{d-s\gamma}}{{d-s\gamma}}  \\
&\leq&\dfrac{2^{d-1}C}{s\gamma-d}|\nn|^{d-s\gamma}
=C_2|\nn|^{d-s\gamma}.
\end{eqnarray*}
Applying the Corollary \ref{coronucleo} we have
\begin{eqnarray*}
\sup_{f\in K*U_p}|| f-sk_{\nn}(f,\cdot)||_q
&\leq& C_3\left( \sum_{|\nl|\geq|\nn|} (a_{\nl})^{s} \right)^{r}
\leq C_{4}|\nn|^{-\gamma+d(p^{-1}-q^{-1})},
\end{eqnarray*}
for $1\leq p\leq 2\leq q\leq \infty$, with $p^{-1}-q^{-1}\geq 2^{-1}$.
\end{prova}

\end{document}